\documentclass[12pt]{article}

\catcode`\@=11

\thicklines
\newskip\Einheit \Einheit=.6cm
\newcount\xcoord \newcount\ycoord
\newdimen\xdim \newdimen\ydim \newdimen\PfadD@cke \newdimen\Pfadd@cke
\def\PfadDicke#1{\PfadD@cke#1 \divide\PfadD@cke by2 
\Pfadd@cke\PfadD@cke \multiply\PfadD@cke by2}
\long\def\LOOP#1\REPEAT{\def\BODY{#1}\ITERATE}
\def\ITERATE{\BODY \let\next\ITERATE \else\let\next\relax\fi \next}
\let\REPEAT=\fi
\def\Punkt{\hbox{\raise-2pt\hbox to0pt{\hss\scriptsize$\bullet$\hss}}}

\def\DuennPunkt(#1,#2){\unskip
  \raise#2 \Einheit\hbox to0pt{\hskip#1 \Einheit
          \raise-1.5pt\hbox to0pt{\hss\tiny$\bullet$\hss}\hss}}
                  
\def\NormalPunkt(#1,#2){\unskip
  \raise#2 \Einheit\hbox to0pt{\hskip#1 \Einheit
          \raise-3pt\hbox to0pt{\hss\large$\bullet$\hss}\hss}}
\def\DickPunkt(#1,#2){\unskip
  \raise#2 \Einheit\hbox to0pt{\hskip#1 \Einheit
          \raise-4pt\hbox to0pt{\hss\Large$\bullet$\hss}\hss}}
\def\Kreis(#1,#2){\unskip
  \raise#2 \Einheit\hbox to0pt{\hskip#1 \Einheit
          \raise-4pt\hbox to0pt{\hss\Large$\circ$\hss}\hss}}
\def\Diagonale(#1,#2)#3{\unskip\leavevmode
  \xcoord#1\relax \ycoord#2\relax
      \raise\ycoord \Einheit\hbox to0pt{\hskip\xcoord \Einheit
         \unitlength\Einheit
         \line(1,1){#3}\hss}}
\def\AntiDiagonale(#1,#2)#3{\unskip\leavevmode
  \xcoord#1\relax \ycoord#2\relax \advance\xcoord by -0.05\relax
      \raise\ycoord \Einheit\hbox to0pt{\hskip\xcoord \Einheit
         \unitlength\Einheit
         \line(1,-1){#3}\hss}}
\def\Pfad(#1,#2),#3\endPfad{\unskip\leavevmode
  \xcoord#1 \ycoord#2 \thicklines\ZeichnePfad#3\endPfad\thinlines}
\def\ZeichnePfad#1{\ifx#1\endPfad\let\next\relax
  \else\let\next\ZeichnePfad
    \ifnum#1=1
      \raise\ycoord \Einheit\hbox to0pt{\hskip\xcoord \Einheit
         \vrule height\Pfadd@cke width1 \Einheit depth\Pfadd@cke\hss}%
      \advance\xcoord by 1
     \else\ifnum#1=2
      \raise\ycoord \Einheit\hbox to0pt{\hskip\xcoord \Einheit
         \unitlength\Einheit
         \line(0,1){1}\hss}
      \advance\xcoord by 0
      \advance\ycoord by 1
 \else\ifnum#1=3
      \raise\ycoord \Einheit\hbox to0pt{\hskip\xcoord \Einheit
         \unitlength\Einheit
         \line(1,1){1}\hss}
      \advance\xcoord by 1
      \advance\ycoord by 1
    \else\ifnum#1=4
      \raise\ycoord \Einheit\hbox to0pt{\hskip\xcoord \Einheit
         \unitlength\Einheit
         \line(1,-1){1}\hss}
      \advance\xcoord by 1
      \advance\ycoord by -1
   \else\ifnum#1=5
      \raise\ycoord \Einheit\hbox to0pt{\hskip\xcoord \Einheit
         \unitlength\Einheit
         \line(2,1){2}\hss}
      \advance\xcoord by 2
      \advance\ycoord by 1
          \else\ifnum#1=6
      \raise\ycoord \Einheit\hbox to0pt{\hskip\xcoord \Einheit
         \unitlength\Einheit
         \line(2,-1){2}\hss}
      \advance\xcoord by 2
      \advance\ycoord by -1
          \else\ifnum#1=7
      \raise\ycoord \Einheit\hbox to0pt{\hskip\xcoord \Einheit
         \unitlength\Einheit
         \line(3,1){3}\hss}
      \advance\xcoord by 3
      \advance\ycoord by 1
          \else\ifnum#1=8
      \raise\ycoord \Einheit\hbox to0pt{\hskip\xcoord \Einheit
         \unitlength\Einheit
         \line(3,-1){3}\hss}
      \advance\xcoord by 3
      \advance\ycoord by -1
    \fi\fi\fi\fi\fi\fi\fi\fi
  \fi\next}
\def\hSSchritt{\leavevmode\raise-.4pt\hbox 
to0pt{\hss.\hss}\hskip.2\Einheit
  \raise-.4pt\hbox to0pt{\hss.\hss}\hskip.2\Einheit
  \raise-.4pt\hbox to0pt{\hss.\hss}\hskip.2\Einheit
  \raise-.4pt\hbox to0pt{\hss.\hss}\hskip.2\Einheit
  \raise-.4pt\hbox to0pt{\hss.\hss}\hskip.2\Einheit}
\def\vSSchritt{\vbox{\baselineskip.2\Einheit\lineskiplimit0pt
\hbox{.}\hbox{.}\hbox{.}\hbox{.}\hbox{.}}}
\def\DSSchritt{\leavevmode\raise-.4pt\hbox to0pt{%
  \hbox to0pt{\hss.\hss}\hskip.2\Einheit
  \raise.2\Einheit\hbox to0pt{\hss.\hss}\hskip.2\Einheit
  \raise.4\Einheit\hbox to0pt{\hss.\hss}\hskip.2\Einheit
  \raise.6\Einheit\hbox to0pt{\hss.\hss}\hskip.2\Einheit
  \raise.8\Einheit\hbox to0pt{\hss.\hss}\hss}}
\def\dSSchritt{\leavevmode\raise-.4pt\hbox to0pt{%
  \hbox to0pt{\hss.\hss}\hskip.2\Einheit
  \raise-.2\Einheit\hbox to0pt{\hss.\hss}\hskip.2\Einheit
  \raise-.4\Einheit\hbox to0pt{\hss.\hss}\hskip.2\Einheit
  \raise-.6\Einheit\hbox to0pt{\hss.\hss}\hskip.2\Einheit
  \raise-.8\Einheit\hbox to0pt{\hss.\hss}\hss}}
\def\SPfad(#1,#2),#3\endSPfad{\unskip\leavevmode
  \xcoord#1 \ycoord#2 \ZeichneSPfad#3\endSPfad}
\def\ZeichneSPfad#1{\ifx#1\endSPfad\let\next\relax
  \else\let\next\ZeichneSPfad
    \ifnum#1=1
      \raise\ycoord \Einheit\hbox to0pt{\hskip\xcoord \Einheit
         \hSSchritt\hss}%
      \advance\xcoord by 1
    \else\ifnum#1=2
      \raise\ycoord \Einheit\hbox to0pt{\hskip\xcoord \Einheit
        \hbox{\hskip-2pt \vSSchritt}\hss}%
      \advance\ycoord by 1
    \else\ifnum#1=3
      \raise\ycoord \Einheit\hbox to0pt{\hskip\xcoord \Einheit
         \DSSchritt\hss}
      \advance\xcoord by 1
      \advance\ycoord by 1
    \else\ifnum#1=4
      \raise\ycoord \Einheit\hbox to0pt{\hskip\xcoord \Einheit
         \dSSchritt\hss}
      \advance\xcoord by 1
      \advance\ycoord by -1
    \fi\fi\fi\fi
  \fi\next}
\def\Koordinatenachsen(#1,#2){\unskip
 \hbox to0pt{\hskip-.5pt\vrule height#2 \Einheit width.5pt depth1 
\Einheit}%
 \hbox to0pt{\hskip-1 \Einheit \xcoord#1 \advance\xcoord by1
    \vrule height0.25pt width\xcoord \Einheit depth0.25pt\hss}}
\def\Koordinatenachsen(#1,#2)(#3,#4){\unskip
 \hbox to0pt{\hskip-.5pt \ycoord-#4 \advance\ycoord by1
    \vrule height#2 \Einheit width.5pt depth\ycoord \Einheit}%
 \hbox to0pt{\hskip-1 \Einheit \hskip#3\Einheit 
    \xcoord#1 \advance\xcoord by1 \advance\xcoord by-#3 
    \vrule height0.25pt width\xcoord \Einheit depth0.25pt\hss}}
\def\Gitter(#1,#2){\unskip \xcoord0 \ycoord0 \leavevmode
  \LOOP\ifnum\ycoord<#2
    \loop\ifnum\xcoord<#1
      \raise\ycoord \Einheit\hbox to0pt{\hskip\xcoord 
\Einheit\Punkt\hss}%
      \advance\xcoord by1
    \repeat
    \xcoord0
    \advance\ycoord by1
  \REPEAT}
\def\Gitter(#1,#2)(#3,#4){\unskip \xcoord#3 \ycoord#4 \leavevmode
  \LOOP\ifnum\ycoord<#2
    \loop\ifnum\xcoord<#1
      \raise\ycoord \Einheit\hbox to0pt{\hskip\xcoord 
\Einheit\Punkt\hss}%
      \advance\xcoord by1
    \repeat
    \xcoord#3
    \advance\ycoord by1
  \REPEAT}
\def\Label#1#2(#3,#4){\unskip \xdim#3 \Einheit \ydim#4 \Einheit
  \def\lo{\advance\xdim by-.5 \Einheit \advance\ydim by.5 \Einheit}%
  \def\llo{\advance\xdim by-.25cm \advance\ydim by.5 \Einheit}%
  \def\loo{\advance\xdim by-.5 \Einheit \advance\ydim by.25cm}%
  \def\o{\advance\ydim by.25cm}%
  \def\ro{\advance\xdim by.5 \Einheit \advance\ydim by.5 \Einheit}%
  \def\rro{\advance\xdim by.25cm \advance\ydim by.5 \Einheit}%
  \def\roo{\advance\xdim by.5 \Einheit \advance\ydim by.25cm}%
  \def\l{\advance\xdim by-.30cm}%
  \def\r{\advance\xdim by.30cm}%
  \def\lu{\advance\xdim by-.5 \Einheit \advance\ydim by-.6 \Einheit}%
  \def\llu{\advance\xdim by-.25cm \advance\ydim by-.6 \Einheit}%
  \def\luu{\advance\xdim by-.5 \Einheit \advance\ydim by-.30cm}%
  \def\u{\advance\ydim by-.30cm}%
  \def\ru{\advance\xdim by.5 \Einheit \advance\ydim by-.6 \Einheit}%
  \def\rru{\advance\xdim by.25cm \advance\ydim by-.6 \Einheit}%
  \def\ruu{\advance\xdim by.5 \Einheit \advance\ydim by-.30cm}%
  #1\raise\ydim\hbox to0pt{\hskip\xdim
     \vbox to0pt{\vss\hbox to0pt{\hss$#2$\hss}\vss}\hss}%
}
\catcode`\@=12

\usepackage{amsmath,amsthm}
\usepackage{amssymb}
\usepackage{pstricks,pst-node}
\usepackage{color}
\usepackage{xspace}
\usepackage{psfig}
\usepackage{epsfig}
\usepackage{graphicx}
\usepackage[colorlinks=true,
linkcolor=green,
filecolor=brown,
citecolor=green]{hyperref}

\definecolor{gray}{rgb}{.221,.221,.221}
\def\gray{\textcolor{gray} }

\def\blue{\textcolor{blue} }

\def\d{\:\vert\:}
\def\a{\ensuremath{\mathcal A}\xspace}

\def\t{\ensuremath{\tau}\xspace}
\def\sc{self-complementary\xspace}
\def\p{\textrm{plane polygon pattern}\xspace}

\def\FPT{\textrm{FPT}}
\newpsobject{showgrid}{psgrid}{subgriddiv=1,griddots=10,gridlabels=6pt}

\begin{document}
\newtheorem{lemma}{Lemma}
\newtheorem{theorem}{Theorem}
\newtheorem{prop}{Proposition}
\newtheorem{cor}{Corollary}
\begin{center}
   
    {\Large Noncrossing Partitions Under Rotation and Reflection      \\ 
}
    \begin{tabbing}
\vspace{10mm}
\hspace*{10mm} \= DAVID CALLAN \hspace*{40mm}   \=  LEN SMILEY \\[-1ex]
\>Department of Statistics \>  Department of Mathematical Sciences \\ [-1ex]
\>University of Wisconsin-Madison  \> University of Alaska Anchorage\\[-1ex]
\>1300 University Ave \> 3211 Providence Drive \\[-1ex]
\>Madison, WI \ 53706-1532 \> Anchorage AK 99508 \\
\>{\bf callan(at)stat.wisc.edu} \> {\bf smiley(at)math.uaa.alaska.edu} \\
\vspace{5mm}
\end{tabbing}

October 30, 2005
\end{center}

\vspace{5mm}

\begin{abstract}
We consider noncrossing partitions of $[n]$ under the action of (i) the
reflection group (of order 2), (ii) the rotation group (cyclic of order n)
and (iii) the rotation/reflection group (dihedral of order 2n). First, we
exhibit a bijection from rotation classes to bicolored plane trees on $n$
edges, and consider its implications. Then we count noncrossing partitions
of $[n]$ invariant under reflection and show that, somewhat surprisingly,
they are equinumerous with rotation classes invariant under reflection.
The proof uses a pretty involution originating in work of Germain
Kreweras.  We conjecture that the ``equinumerous" result also holds for
arbitrary partitions of [$n$].
\end{abstract}

\vspace*{10mm}
{\large \textbf{1 \quad Introduction}  }

A \emph{noncrossing} partition of $[n]=\{1,2,\ldots,n\}$ is one for which 
no quadruple $a<b<c<d$ has $a,c$ in one block and $b,d$ in 
another. This implies that if the elements of $[n]$ are situated around a circle with $1\to 2\to 3\to\dots\to n\to 1$ forming a
 cycle, and neighboring elements 
within each block are joined by line segments, then no line segments 
cross one another (see figure below).

The map $i \mapsto i+1\,$(mod\:$n$) on $[n]$ induces a map---the 
rotation operator $R$---on partitions $\pi$ of $[n]$. Equivalence under 
repeated application of $R$ divides them into rotation 
classes: $A(\pi)=\{R^{i}(\pi)\}_{i\ge 1}$. 
The \emph{complement} of a partition $\pi$ of $[n]$ is 
$C(\pi):=n+1-\pi$ (elementwise).
It is easy to check that $C\circ R =R^{-1}\circ C$ and so the 
complement operation permutes rotation classes. We say a partition 
$\pi$ is 
\emph{self-complementary} if $C(\pi)=\pi$ and a rotation class $A$ 
is \sc if $C(A)=A$. As we will see, a \sc rotation class need not 
contain any \sc partitions.

The operations rotation and complementation both preserve the 
noncrossing (NC) property of partitions. In particular, a rotation class 
consists entirely of NC partitions if it contains a single one. A NC 
rotation class may be represented by a polygon diagram with the labels 
removed, we'll call it an \emph{NC Polygon (or Partition) Pattern} (NCPP).
\begin{center}

\begin{pspicture}(-8,0)(8,4)
    
\pscircle(-4,2){2}
\psdots(-2,2)(-3,3.73)(-3,.27)(-5,.27)(-5,3.73)(-6,2)
\pspolygon(-3,3.73)(-3,.27)(-5,.27)
\psline(-5,3.73)(-6,2)
\rput(-2,2){$\bullet$}

\pscircle(4,2){2}
\psdots(2,2)(3,3.73)(3,.27)(5,.27)(5,3.73)(6,2)
\pspolygon(5,3.73)(3,.27)(5,.27)
\psline(3,3.73)(2,2)
\rput(6,2){$\bullet$}

\rput(-2.9,3.95){\textrm{{\footnotesize 1}}}
\rput(-1.8,2){\textrm{{\footnotesize 2}}}
\rput(-2.9,.05){\textrm{{\footnotesize 3}}}
\rput(-5.1,.05){\textrm{{\footnotesize 4}}}
\rput(-6.2,2){\textrm{{\footnotesize 5}}}
\rput(-5.1,3.95){\textrm{{\footnotesize 6}}}

\end{pspicture}

\end{center} 
\begin{center}
    {\footnotesize polygon diagram of the NC partition\hspace*{21mm} polygon 
    diagram of its rotation class}\\ \vspace*{-1mm}
    {\footnotesize  \hspace*{4mm}134-2-56: labels fixed in place 
    \hspace*{28mm} is an NCPP: no labels, rotate at will}
\end{center}

Clearly, a NC partition is \sc if its labeled polygon diagram is 
invariant when flipped across a vertical line, and a NC rotation 
class is \sc if its \p is achiral, that is,  
invariant when flipped over (across any line).

\begin{figure}[h]
\includegraphics[scale=0.6, trim= 0 300 0 300]{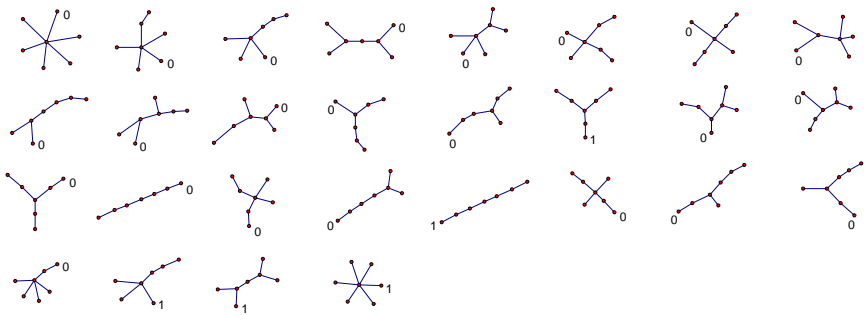}
\caption{Bicolored 6-edge plane trees}
\label{F66}
\end{figure}

A \emph{bicolored plane tree} is a plane tree (no root, no labels) in 
which each vertex is colored white or yellow (say) in such a way that adjacent 
vertices get different colors. The color of one vertex determines 
that of all the others (see Figure \ref{F66}). 
We exhibit a bijection from NC rotation classes ($n$ points) to bicolored plane 
trees ($n$ edges) in \S 2. We count \sc NC partitions in \S 3 and show they are 
equinumerous with \sc NC rotation classes, equivalently, achiral 
NC polygon patterns in \S 4. Remarks, figures, and a conjecture comprise \S 5.
The Appendix contains enumerations and tables.
The enumeration of bicolored plane trees
(\htmladdnormallink{A054357}{http://www.research.att.com:80/cgi-bin/access.cgi/as/njas/sequences/eisA.cgi?Anum=A054357})
has been significantly generalized in \cite{BousquetCacti}. 
We note that (scaled) 
Shabat polynomials \cite{Adrianov,Bordeaux} are also counted by
\htmladdnormallink{A054357}{http://www.research.att.com:80/cgi-bin/access.cgi/as/njas/sequences/eisA.cgi?Anum=A054357}.

\vfill
\eject

{\large \textbf{2 \quad Bijection}  }

The bijection from NC polygon diagrams to bicolored binary trees 
is depicted in the Figure below. In forming the polygon diagram of 
a NC partition, the convex hull of each block of size $k$ forms a
$k$-sided yellow polygon even for $k=1,2$ by liberal interpretation of
``polygon'' as illustrated. The polygons are disjoint because the partition
is noncrossing.  Ignoring the labels and considering the configuration of
polygons only up to rotation, it represents an NC polygon pattern (NCPP). 

\begin{center}

\begin{pspicture}(-6,-6)(6,6)
\psset{unit=1cm} 
 \SpecialCoor
\pswedge[linecolor=lightgray,fillstyle=solid,fillcolor=yellow](5;288){.4}{25}{190} 
\pswedge[linecolor=lightgray,fillstyle=solid,fillcolor=yellow](5;180){.4}{280}{80}

\pscircle[linecolor=red](0,0){5}

\pspolygon[linecolor=lightgray,fillstyle=solid,fillcolor=yellow](5;72)(5;108)(5;252)
\pspolygon[linecolor=lightgray,fillstyle=solid,fillcolor=yellow](5;216)(4.25;180)(5;144)(3.95;180)
\pspolygon[linecolor=lightgray,fillstyle=solid,fillcolor=yellow](5;0)(5;36)(5;324)
\psline(4.8;18)(4.3;0)(4.8;342)
\psline(4.8;90)(1.5;120)(2;0)(4.3;0)
\psline(1.5;120)(2.5;180)(4.1;180)(4.5;185)(4.85;180)
\psline(2;0)(4.85;288)

\gray{
\rput(5.3;0){\textrm{{\footnotesize 3}}}
\rput(5.3;36){\textrm{{\footnotesize 2}}}
\rput(5.3;72){\textrm{{\footnotesize 1}}}
\rput(5.3;108){\textrm{{\footnotesize 10}}}
\rput(5.3;144){\textrm{{\footnotesize 9}}}
\rput(5.3;180){\textrm{{\footnotesize 8}}}
\rput(5.3;216){\textrm{{\footnotesize 7}}}
\rput(5.3;252){\textrm{{\footnotesize 6}}}
\rput(5.3;288){\textrm{{\footnotesize 5}}}
\rput(5.3;324){\textrm{{\footnotesize 4}}}
}

\psdots(2;0)(4.3;0)
(2.5;180)(4.1;180)(4.5;185)(4.85;180)
(4.85;90)(4.85;18)(4.85;342)(1.5;120)(4.85;288)

\end{pspicture}
\end{center}

\vspace*{-15mm}

\begin{center}
    {\footnotesize An NC polygon diagram with superimposed bicolored plane tree}\\ \vspace*{-1mm}
    {\footnotesize  polygon sides $\leftrightarrow$ tree edges }
\end{center}

The bijection is clear: place a vertex in each region of the circle, both
yellow and white.  Join vertices in adjacent regions by edges.  Then allow
each vertex to inherit the color of the region it's in to get the desired
bicolored plane tree.  As for invertibility, the ``star'' tree with yellow leaves joined to a 
white center corresponds to the NCPP all of whose polygons are of the 
degenerate one-sided type. Otherwise the tree has a non-leaf yellow 
vertex and its inverse is formed recursively as illustrated below.
\begin{center}

\begin{pspicture}(-4,-4)(4,3)
\psset{unit=.5cm}   
\pscircle[linecolor=red](0,0){5.66}
\pscircle[fillstyle=solid,fillcolor=yellow](0,0){0.3}
\pspolygon[linewidth=0.5mm,linecolor=red](-4,-4)(-4,4)(4,4)(4,-4)
\psline[linewidth=0.3mm](-4.5,0)(4.5,0)
\psline(0,-4.5)(0,4.5)

\rput(0,5){\textrm{{\footnotesize $T_{1}$}}}
\rput(5,0){\textrm{{\footnotesize $T_{2}$}}}
\rput(0,-5){\textrm{{\footnotesize $T_{3}$}}}
\rput(-5,0){\textrm{{\footnotesize $T_{4}$}}}
\end{pspicture}

\vspace*{-15mm}

\end{center} 
\begin{center}
    {\footnotesize the polygon surrounding an internal yellow vertex; recursively}\\ \vspace*{-1mm}
    {\footnotesize  construct polygons for the trees $T_{i}$ on 
    their corresponding arcs}
\end{center}
There are two immediate consequences of this bijection. The first one 
is a special case of Theorem 10 in \cite{BousquetCacti}.
\begin{itemize}
    \item  The Catalan numbers
    (\htmladdnormallink{A000108}{http://www.research.att.com:80/cgi-bin/access.cgi/as/njas/sequences/eisA.cgi?Anum=A000108})
    count bicolored plane trees with a 
    distinguished edge. This is because the positioning of the labels 
    on the circle can be captured by associating an edge in the tree 
    with label 1, say the first tree edge encountered travelling 
    clockwise from 1 around the polygon incident with 1. Thus the 
    distinguished-edge bicolored plane trees on $n$ edges are in correspondence 
    with ordinary NC partitions of $[n]$, counted by $C_{n}$ \cite{ec2}.
    \item  On NC partitions, as well as on circular NC partitions, the statistics ``\#\,singletons'' 
    and ``\#\,adjacencies'' have the same distribution, in fact a symmetric joint distribution.
    This is due to the correspondences yellow leaf $\leftrightarrow$ 
    singleton block, and white leaf $\leftrightarrow$ adjacency, that is, 
    two consecutive 
    elements of $[n]$ in the same block. (Of course consecutive is taken 
    here in the circular sense, so $n$ and 1 are considered consecutive.) The symmetry 
    of the joint distribution holds for unrestricted partitions too \cite{conjpart}.
   
\end{itemize} 

\vfill
\eject

{\large \textbf{3 \quad Counting Self-Complementary NC Partitions}  }

Using paths from $(0,0)$ with Upsteps $(1,1)$ and Downsteps $(1,-1)$ (see, e.g., \cite{Deutsch}), 
we will show that the number of \sc NC partitions of $[n]$ is $\binom{n}{\lfloor n/2 
\rfloor}$.
There is a bijection from Dyck 
$n$-paths (Quadrant 1, ending at $(2n,0)$) to NC partitions of $[n]$ that sends \#\,peaks to 
\#\,blocks. Given a Dyck $n$-path, number its upsteps left to right and
then give each downstep the number of its matching upstep. The numbers 
on each descent (maximal sequence of contiguous downsteps) form the 
blocks of the corresponding NC partition.

\vbox{
$$
\Einheit=0.6cm.
\Pfad(-12,0),343334344333443444334344\endPfad
\SPfad(-12,0),111111111111111111111111\endSPfad
\blue{
\Label\o{ \textrm{{\scriptsize 1}}    }(-10.3,0.3)
\Label\o{ \textrm{{\scriptsize 4}}    }(-6.3,2.3)
\Label\o{ \textrm{{\scriptsize 5}}    }(-4.3,2.3)
\Label\o{ \textrm{{\scriptsize 3}}    }(-3.3,1.3)
\Label\o{ \textrm{{\scriptsize 8}}    }(0.8,3.3)
\Label\o{ \textrm{{\scriptsize 7}}    }(1.8,2.3)
\Label\o{ \textrm{{\scriptsize 9}}    }(3.8,2.3)
\Label\o{ \textrm{{\scriptsize 6}}    }(4.8,1.3)
\Label\o{ \textrm{{\scriptsize 2}}    }(5.8,0.3)
\Label\o{ \textrm{{\scriptsize 11}}    }(8.8,1.3)
\Label\o{ \textrm{{\scriptsize 12}}    }(10.8,1.3)
\Label\o{ \textrm{{\scriptsize 10}}    }(11.8,0.3)
}
\Label\o{ \textrm{{\tiny 1}}    }(-11.7,0.3)
\Label\o{ \textrm{{\tiny 2}}    }(-9.7,0.3)
\Label\o{ \textrm{{\tiny 3}}    }(-8.7,1.3)
\Label\o{ \textrm{{\tiny 4}}    }(-7.7,2.3)
\Label\o{ \textrm{{\tiny 5}}    }(-5.7,2.3)
\Label\o{ \textrm{{\tiny 6}}    }(-2.7,1.3)
\Label\o{ \textrm{{\tiny 7}}    }(-1.7,2.3)
\Label\o{ \textrm{{\tiny 8}}    }(-0.7,3.3)
\Label\o{ \textrm{{\tiny 9}}    }(2.3,2.3)
\Label\o{ \textrm{{\tiny 10}}    }(6.3,0.3)
\Label\o{ \textrm{{\tiny 11}}    }(7.3,1.3)
\Label\o{ \textrm{{\tiny 12}}    }(9.3,1.3)
\Label\u{ \textrm{ {\footnotesize  number downsteps,}}  }(0,-0.2)
\Label\u{ \textrm{ {\footnotesize each with the number of its matching upstep}} }(0,-0.8)
\DuennPunkt(-12,0)
\DuennPunkt(-11,1)
\DuennPunkt(-10,0)
\DuennPunkt(-9,1)
\DuennPunkt(-8,2)
\DuennPunkt(-7,3)
\DuennPunkt(-6,2)
\DuennPunkt(-5,3)
\DuennPunkt(-4,2)
\DuennPunkt(-3,1)
\DuennPunkt(-2,2)
\DuennPunkt(-1,3)
\DuennPunkt(0,4)
\DuennPunkt(1,3)
\DuennPunkt(2,2)
\DuennPunkt(3,3)
\DuennPunkt(4,2)
\DuennPunkt(5,1)
\DuennPunkt(6,0)
\DuennPunkt(7,1)
\DuennPunkt(8,2)
\DuennPunkt(9,1)
\DuennPunkt(10,2)
\DuennPunkt(11,1)
\DuennPunkt(12,0)
$$
}
Partition downstep labels by descents to get 
\[
\blue{1\,-\,4\,-\,5\,3\,-\,8\,7\,-\,9\,6\,2\,-\,11\,-\,12\ 10},
\]
a noncrossing partition with arc diagram 
\begin{center}

\begin{pspicture}(-8,0)(8,1)

\dotnode(-6,0){A}
\dotnode(-5,0){B}
\dotnode(-4,0){C}
\dotnode(-3,0){d}
\dotnode(-2,0){e}
\dotnode(-1,0){f}
\dotnode(0,0){g}
\dotnode(1,0){h}
\dotnode(2,0){j}
\dotnode(3,0){k}
\dotnode(4,0){m}
\dotnode(5,0){n}

\rput(-6,-0.3){{\footnotesize 1}}
\rput(-5,-0.3){{\footnotesize 2}}
\rput(-4,-0.3){{\footnotesize 3}}
\rput(-3,-0.3){{\footnotesize 4}}
\rput(-2,-0.3){{\footnotesize 5}}
\rput(-1,-0.3){{\footnotesize 6}}
\rput(0,-0.3){{\footnotesize 7}}
\rput(1,-0.3){{\footnotesize 8}}
\rput(2,-0.3){{\footnotesize 9}}
\rput(3,-0.3){{\footnotesize 10}}
\rput(4,-0.3){{\footnotesize 11}}
\rput(5,-0.3){{\footnotesize 12}}

\pscurve(-5,0)(-4,0.8)(-2,0.8)(-1,0)
\pscurve(-4,0)(-3.7,.3)(-2.3,.3)(-2,0)
\pscurve(-1,0)(-0.5,.5)(1.5,.5)(2,0)
\pscurve(3,0)(3.3,.3)(4.7,.3)(5,0)
\pscurve(0,0)(0.3,.2)(0.7,.2)(1,0)
\end{pspicture}

\end{center}
Under this bijection, peak downsteps correspond to largest block 
elements, and downsteps returning the path to ground level correspond to 
smallest elements in maximal blocks (a maximal block is one whose 
arcs would get wet if it rained, here there are 3 such: 1, 2\:6\:9, 10\:12).
\begin{theorem} The number of \sc NC partitions of $[n]$ is
$ \vert \a_{n} \vert =\binom{n}{\lfloor n/2\rfloor}.$
\end{theorem}
\noindent \textbf{Proof}\quad We give a bijective proof for $n$ even. (The 
case $n$ odd is similar and is omitted.) So suppose $n=2m$. The right 
hand side clearly counts paths of $m$ upsteps and $m$ downsteps (\emph{balanced} $m$-paths). 
Now a NC 
partition $\pi$ of $[2m]$ induces a partition $\tau$ of $[m]$ by 
intersecting its blocks with $[m]$. For $\t_{i}$ a block of \t, set 
$\overline{\t_{i}}=2m+1-\t$ (elementwise). If $\pi$ is \sc then so 
is $\t$ and each block of $\pi$ has one of the three forms $\t_{i},\ 
\overline{\t_{i}}$ or $\t_{i}\cup \overline{\t_{i}}$ for some block $
\t_{i}$ of \t. The first two forms come in complementary pairs, the 
last form is permissible only if $\t_{i}$ is a maximal block of \t 
(else $\pi$ would have a crossing). So \sc NC partitions $\pi$ of 
$[2m]$ correspond to NC partitions \t of $[m]$ in which each maximal 
block may (or not) be marked: a mark on $\t_{i}$ indicating that  
$\t_{i}\cup \overline{\t_{i}}$ is a block in $\pi$, the absence of a 
mark indicating that $\t_{i},\ \overline{\t_{i}}$ are separate blocks 
of $\pi$. Using the NC partition $\leftrightarrow$ Dyck path 
correspondence above, these marked objects correspond in turn to Dyck 
$m$-paths with returns (to ground level) available for marking. 
Returns split a Dyck path into its components (Dyck subpaths whose 
only return is at the end). Flip over each component that terminates 
at a marked return to obtain a balanced $m$-path. This is the desired bijection 
from \sc NC partitions of $[2m]$ to balanced $m$-paths. \qed

\vspace*{10mm}
{\large \textbf{4 \quad Counting Achiral NC Polygon Patterns}  }

\begin{theorem} The set $\a_{n}$ of achiral NC Polygon Patterns $($or \sc NC rotation 
classes$)$ of $[n]$ is equinumerous with the set of \sc NC partitions 
of $[n]$, and hence $\vert \a_{n} \vert =\binom{n}{\lfloor n/2\rfloor}$.
\label{achiral}
\end{theorem}
To prove this, recall two related operations on NC partitions 
\cite{krewerasNC,simionNC} defined using polygon diagrams as 
illustrated in Figures \ref{F8},\ref{F9},\ref{F10} below. 


In both cases, new vertices (in blue) interleave the 
old vertices (in black) but their labelings differ. The new labels 
are then formed into maximal blocks subject only to: new polygons are 
disjoint from the old ones.

\begin{figure}[h]
\centering
\begin{minipage}{.25\textwidth}
\centering\includegraphics[scale=0.5, trim= 150 200 100 200]{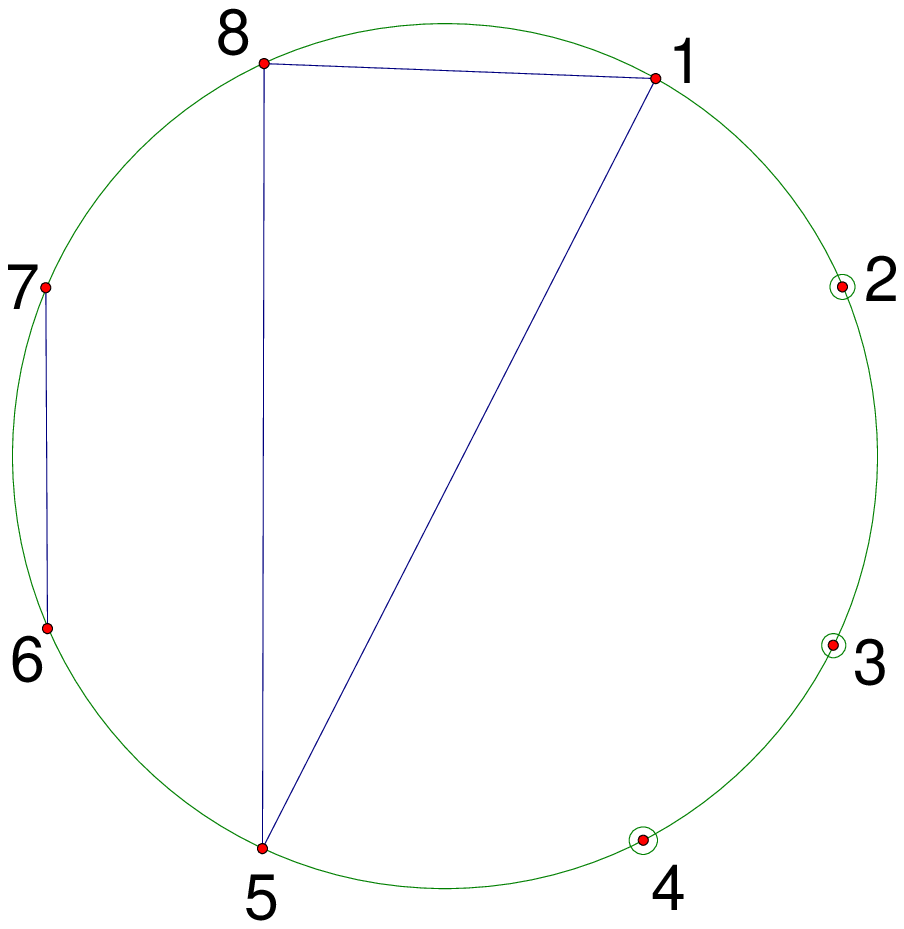}\hfill
\caption{NC partition $\pi$ by polygons}\hfill
\label{F8}
\end{minipage}
\hfill
\begin{minipage}{.25\textwidth}
\centering\includegraphics[scale=0.5, trim= 150 200 100 200]{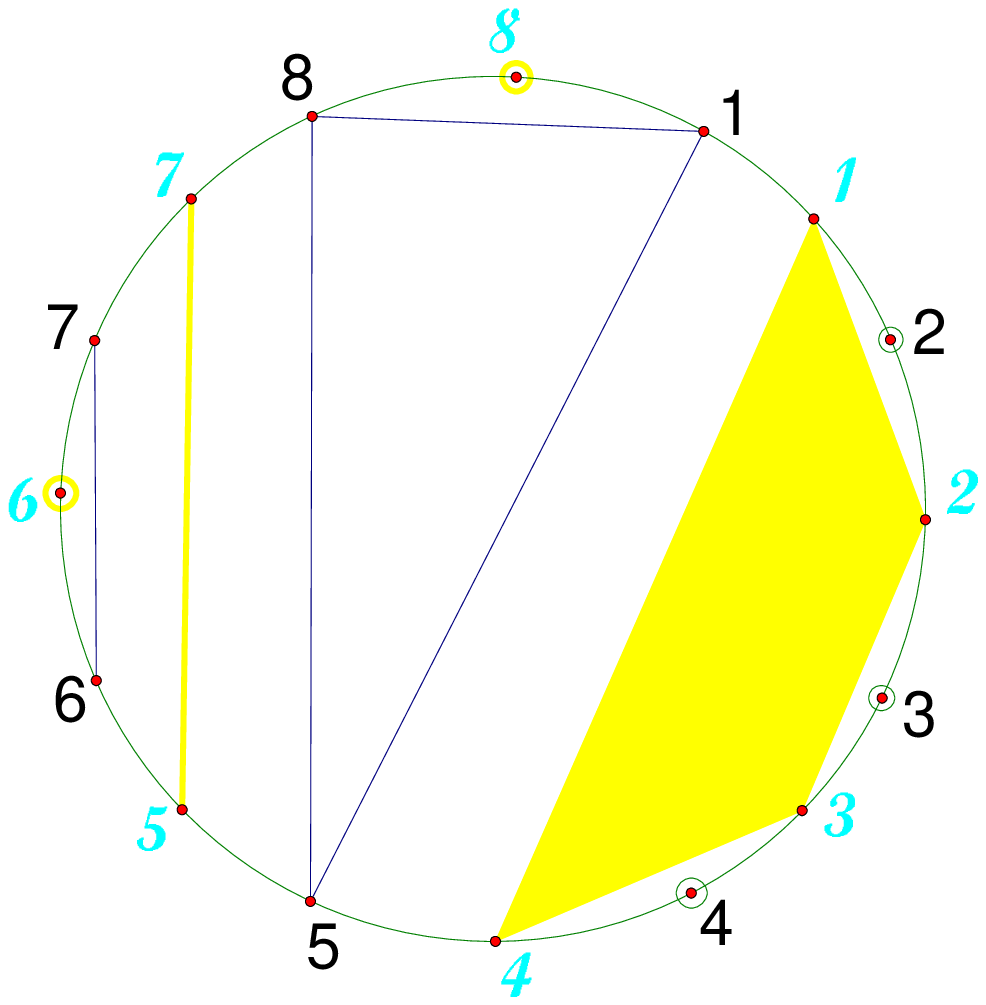}\hfill
\caption{$\pi$ with $H(\pi )$}\hfill
\label{F9}
\end{minipage}
\hfill
\begin{minipage}{.25\textwidth}
\centering\includegraphics[scale=0.5, trim= 150 200 100 200]{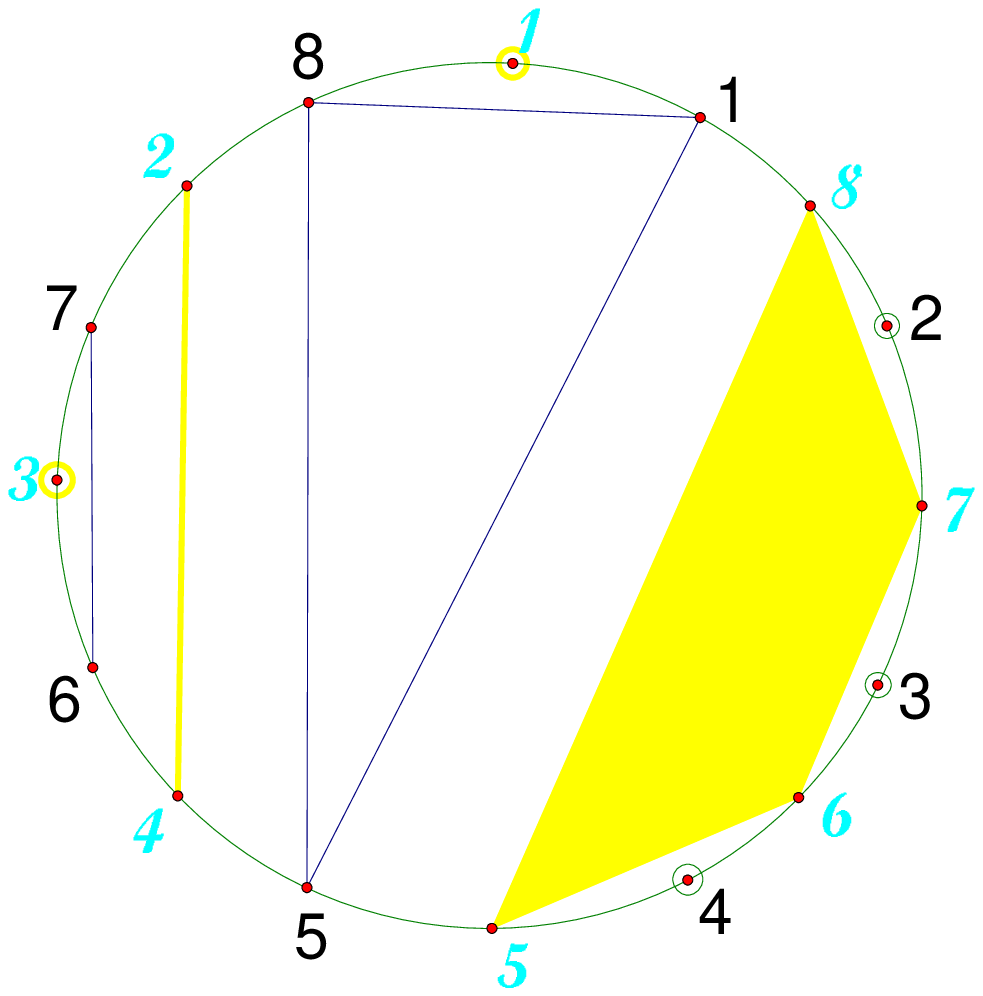}
\caption{$\pi$ with $T(\pi )$}
\label{F10}
\end{minipage}
\hfill
\end{figure}

It is clear from their defining diagrams that 
\[
H^{2}=R^{-1},\quad T^{2}=I,\quad T=CH
\]
and it is not hard to see that $TR=R^{-1}T$ and hence, by induction, 
$TR^{i}=R^{-i}T$ for all $i$. The following result is key to the 
bijection establishing the Theorem.

\begin{prop}
    CT=TRC
\label{keyid}    
\end{prop}
\vspace*{-3mm}
\noindent \textbf{Proof}\quad $H^{2}=R^{-1} \Rightarrow H=R^{-1}H^{-1} \Rightarrow 
CT=R^{-1}TC=TRC$. \qed

For $n\ge 3,\ H \ne R$ and we note in passing that the operations $C,T,H,R$ on NC 
partitions of $[n]$ generate a dihedral group $D_{2n}$ (of $4n$ 
elements) with 
presentation $\langle H,C\: 
:\:H^{2n}=C^{2}=I,\ CHC^{-1}=H^{-1}\rangle$.

Next we define the notion of \emph{complement order} on partitions in 
achiral NC rotation classes. Suppose $A \in \a_{n}$ and $\pi$ is a partition 
in $A$. Then $C(\pi) \in A$ (because $A$ is achiral) and so $C(\pi)=R^{i}(\pi)$ for some $i\ge 
1$ ($i=n$ will do if $C(\pi)=\pi$). Define the complement order of 
$\pi$ to be the minimal such $i\ge 1$.

\begin{lemma} \textrm{ }\\
\hspace*{1mm} $($i\,$)$ An achiral rotation class $A$ contains at most 2 
\sc partitions.\\
\hspace*{1mm}$($ii\,$)$ If $\vert A \vert$ is odd, then $A$ contains exactly one \sc 
partition. \\
$($iii\,$)$ If $\vert A \vert$ is even, then either every partition in $A$ 
has even complement order or every partition in $A$ 
has odd complement order. In the former case, $A$ contains 2 \sc 
partitions; in the latter case, none.
\label{3props}
\end{lemma}
The proof is deferred. Theorem \ref{achiral} will follow from this lemma if we can show that, among 
even-cardinality achiral NC rotation classes $A$ in $\a_{n}$, there are 
just as many associated with even complement order as with odd. (Of 
course, $\vert A \vert$ even implies $n$ even.) We claim the transpose 
$T$ is a bijection, indeed an involution, that interchanges these two
families. To see this, first suppose that $A \in \a_{n}$ has even 
cardinality, say $\vert A 
\vert=2s$, and $\pi \in A$ has even complement order, say
$C(\pi)=R^{2m}\pi$. Then, using Proposition \ref{keyid},
\[
C(T\pi)=TR\,C\pi=TR^{2m+1}\pi=R^{-(2m+1)}T\pi=R^{2s-2m-1}(T\pi)
\]
and $T\pi$ has odd complement order. The other direction is 
similar, the desired bijection is established, and Theorem \ref{achiral} follows.

\noindent \textbf{Proof of Lemma 1}\quad Suppose a 
rotation class $A$ contains a \sc partition $\pi$. Then the complement of 
every other element of $A$ is given by 
\begin{equation}
    CR^{i}\pi = R^{-i}C\pi=R^{-i}\pi
    \label{comp}
\end{equation}
Now suppose $R^{i}\pi\in A$ is also \sc. It follows from (\ref{comp}) 
that $R^{2i}\pi=\pi$. Set $t=\vert A \vert$ so that $R^{j}\pi=\pi 
\Rightarrow t\d j.$ Hence $t\d 2i$.

If $t$ is odd, then $t\d i $ and $R^{i}\pi=\pi$, implying that $\pi$ 
is the only \sc partition in $A$. If $t$ is even, say $t=2s$, then 
$s\d i$ and $R^{i}\pi$ is one of $R^{s}\pi$ and
$R^{2s}\pi=\pi$. These facts establish part (i) and the ``at most one'' half of part 
(ii).

For the ``at least one'' half of part (ii), suppose $\vert A\vert$ is odd. Take $\pi \in A$. Since $A$ is 
achiral, $C\pi=R^{k}\pi$ for some $k$ and so the complement of each 
element of $A$ is given by $CR^{i}\pi=R^{k-i}\pi$. If $k$ is even, 
then $i=k/2$ makes $R^{i}\pi$ \sc. On the other hand, if $k$ is odd, 
say $k=2\ell+1$ and $\vert A \vert=2s+1$, then $i=\ell-s$ makes $R^{i}\pi$ 
\sc. This establishes part (ii). 

For part (iii), let $t:=\vert A \vert$ be even. First, suppose some 
$\pi \in A$ has even complementary order $k$: $C(\pi)=R^{k}\pi$. Then 
$C(R^{i}\pi)=R^{k-i}\pi=R^{k-2i}R^{i}\pi$ and the powers of $R$ that 
fix $R^{i}\pi$ are all $\equiv k-2i\,$(mod\,$t$) and hence even. Thus 
every element of $A$ has even complementary order and $R^{i}\pi$ is 
\sc for $i=k/2$ and $i=(k+t)/2$. Similarly, if some element $\pi$ of $A$ 
has odd complementary order, then they all do, and the equation 
$C(R^{i}\pi)=R^{i}\pi$ has no solution. \qed

\begin{figure}[h]
\includegraphics[scale=0.45, trim= 0 300 0 80]{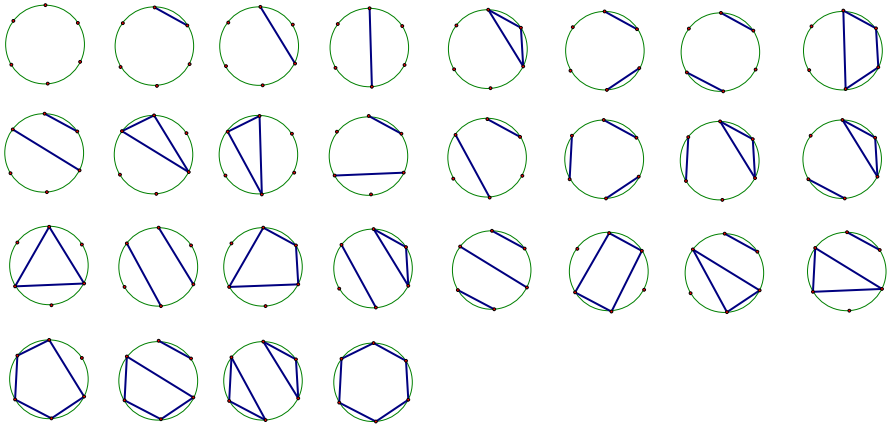}
\caption{NC Patterns, n=6}
\label{F0}
\end{figure}

\vspace*{10mm}
{\large \textbf{5 \quad Concluding Remarks}  }
\vspace*{-1mm}
\begin{enumerate}
    \item  The transpose defined in Figure \ref{F10} above coincides with the restriction 
    to NC partitions
    of the conjugate \cite{conjpart} defined on all partitions of 
    $[n]$. In particular, the algorithmic definition of conjugate 
    given in 
    \cite{conjpart} provides a practical way to compute the transpose.

    \item  An analog of Theorem \ref{achiral} appears to hold for arbitrary 
    partitions: the number of \sc rotation classes of partitions on 
    $[n]$ coincides with the number of \sc partitions of $[n]$. The proof of 
    Lemma \ref{3props} goes through unchanged (it does not use the NC property). 
    Unfortunately, the conjugate does not serve in the role  of 
    transpose to interchange the two relevant families in this 
    larger setting, and it would be 
    interesting to find an extension of the transpose that does.

    \item  Figures \ref{F0} and \ref{F66} show the 28 bijective pairs (by geographic position) 
    of NCPP's 
    and bicolored plane trees for $n=6$. To avoid clutter, only one leaf per tree displays its color.
    Singleton (yellow) leaves are colored ``$0$''; adjacency (white) leaves are colored ``$1$''.

\end{enumerate}

\vfill
\vfill
\eject

\vspace*{10mm}
 
{\large \textbf{Appendix: \quad Enumeration Formula for Unlabeled NC Partitions and Bicolored Plane Trees (almost \emph{ab initio})}  }

The partitions of the set $[n]=\{ 1,2,\dots ,n\}$ (the decompositions of $[n]$ as a union of pairwise-disjoint,
 non-empty subsets) are counted by the sequence of Bell numbers.
If the elements of $[n]$ are regarded as the set of labels of $n$ otherwise indistinguishable
objects, the unlabeled enumeration of partitions of these objects is the same as counting the 
partitions of the integer $n$.

Non-crossing partitions (cf. Introduction) are beautiful, and have been closely studied. 
Motzkin noted \cite[last sentence]{MotzkinBAMS} that the number of labeled non-crossing partitions, as a function of $n$,
 satisfies the Catalan recurrence, and in fact these are counted by the sequence of Catalan numbers.

In the unlabeled case, there are two candidate sequences: the leading one regards the circle as embedded in a plane with the points evenly spaced and counts non-crossing partitions inequivalent under rotations of the circle (these are the \emph{NC Partition Patterns}); the second identifies two partition classes counted in the first which are the same after reflection across a diameter of the circle (we might call these classes \emph{chirally inequivalent NC Partition Patterns}). Motzkin  \cite[penultimate sentence]{MotzkinBAMS} gave the beginning of the latter sequence as $1,2,3,6,9,24$. This contains an (almost certainly clerical) error: the value for $n=5$ should be $10$, not $9$.

Because the ``label/unlabel'' paradigm has been invoked, the sequence $NCPP(n)$ fits into one of at least two competing enumerative analogies: if the circle is discarded entirely, then Bell numbers\ :\ partitions of $n$\ ::\ Catalan numbers\ :\ $NCPP(n)$; if only the NC requirement is relaxed, then  Bell\ :\ possibly crossing partition patterns (\htmladdnormallink{A084423}{http://www.research.att.com:80/cgi-bin/access.cgi/as/njas/sequences/eisA.cgi?Anum=A084423})\ ::\ Catalan\ :\ $NCPP(n)$.

The essential fact needed in the direct enumeration of NC partition patterns was proven by V. Reiner 
\cite{reinerNC}: if $n=kd$, $k\geq 2$,
 and the points are labelled $a_1,a_2,\dots ,a_n$, then
the number of NC partitions having the property ``$a_i$ and $a_j$ are in the same part if and only if 
$a_{i+d\pmod n}$ and $a_{j+d\pmod n}$ are in the same part'' is $\binom{2d}{d}$. We refer to such NC partitions
as \emph{d-clickable} (suggested by analogy to clicking a physical dial with $n$ positions through $d$ positions
and arriving at the same partition). 
Reiner gives two proofs; for the 
convenience of the reader we informally describe the bijection used in one of them. Using $k$ distinct colors for the $a$'s,
relabel the points consecutively to make $k$ monocolor intervals subscripted $1,\dots ,d$. 
Consider the 'unwrapped' doubly-infinite 
sequence
$$\dots\mathbf{a}_{d-1},\mathbf{a}_{d},a_{1},{a}_{2},\dots ,{a}_{d},A_{1},\dots ,{A}_{d},\mathsf{a}_{1},\mathsf{a}_{2}\dots\dots\dots ,\mathbf{a}_{d-1},\mathbf{a}_{d},a_{1},\dots$$
(here font/case is used to denote color). On the circle, find a part of the partition, say of size $p$, consisting entirely of
a consecutive set of points. This will always be possible for an NC partition. 
In a set $L$ place the subscript of the first (clockwise) element in the chosen part,
 and  in a set $R$ the subscript of the last. This will be possible if the part is proper. Remove all elements with subscripts equal to those in this part
 from the 
doubly-infinite sequence and from the circle. The resulting sequence still consists of equal length monocolor intervals
in the $k$ colors, and (after equispacing the remaining points on the circle) is $(d-p)$-clickable. Repeat the process until no such proper part remains to be chosen, at which time the sets $L$ and $R$ are equinumerous, but otherwise arbitrary, subsets of $[d]$. The number of ways of specifying
such an $L$ and $R$ is easily seen to be $\binom{2d}{d}$.

To reverse the process, consider a copy of the original doubly-infinite sequence and for each element of 
 $L$ (resp. $R$) place a Left (resp. Right) parenthesis to the left (resp. right)
of each symbol in the sequence with subscript equal to this element. When $L$ and $R$ are exhausted the partition
may be decoded from the parenthesized string in the usual manner.

\begin{figure}[h]
\begin{minipage}{2.6 in}
\includegraphics[scale=0.35, trim= 0 30 0 80]{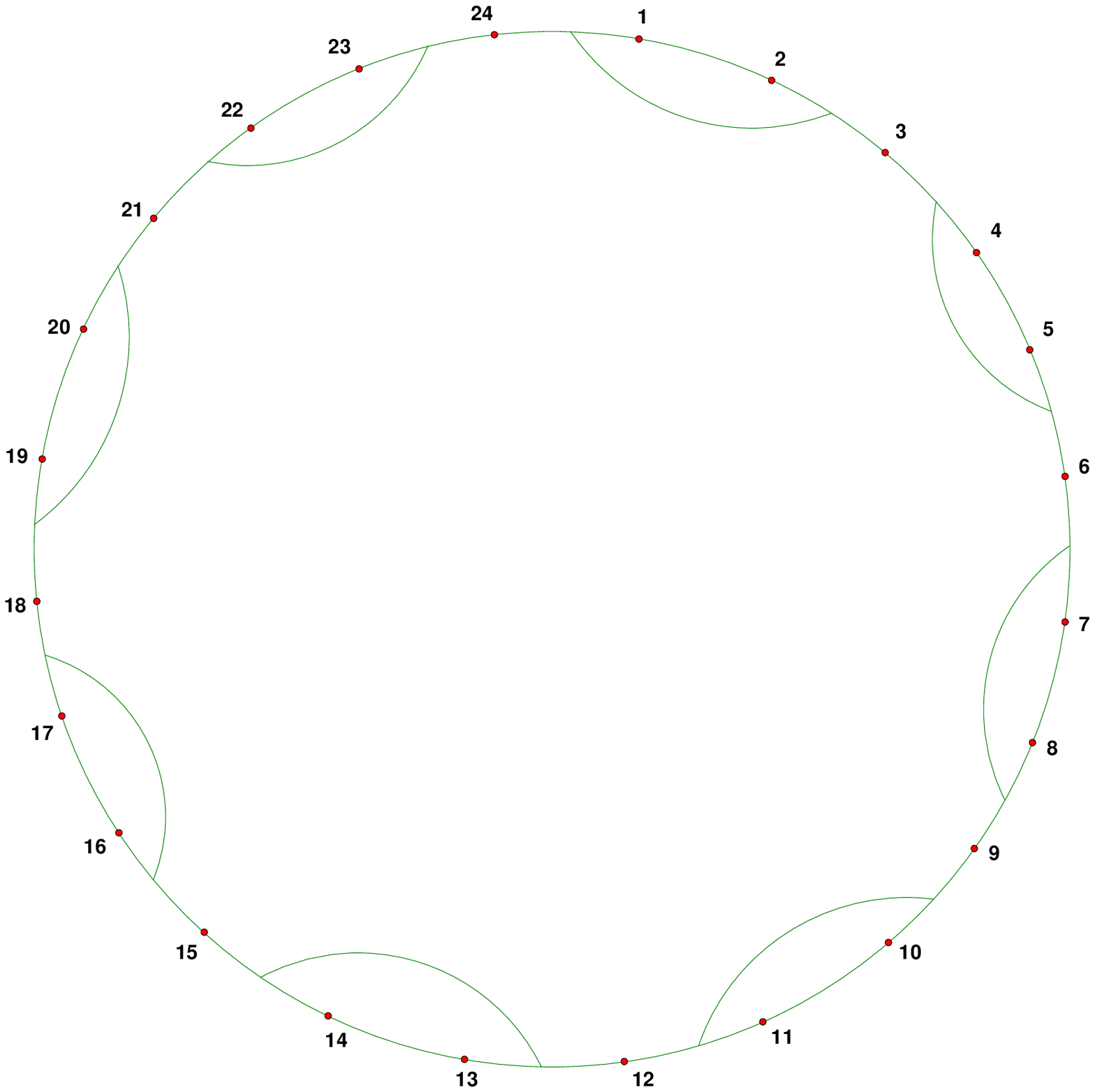}
\caption{A 3-clickable partition}
\label{F3}
\end{minipage}
\hfill
\begin{minipage}{2.6 in}
\includegraphics[scale=0.35, trim= 0 30 0 80]{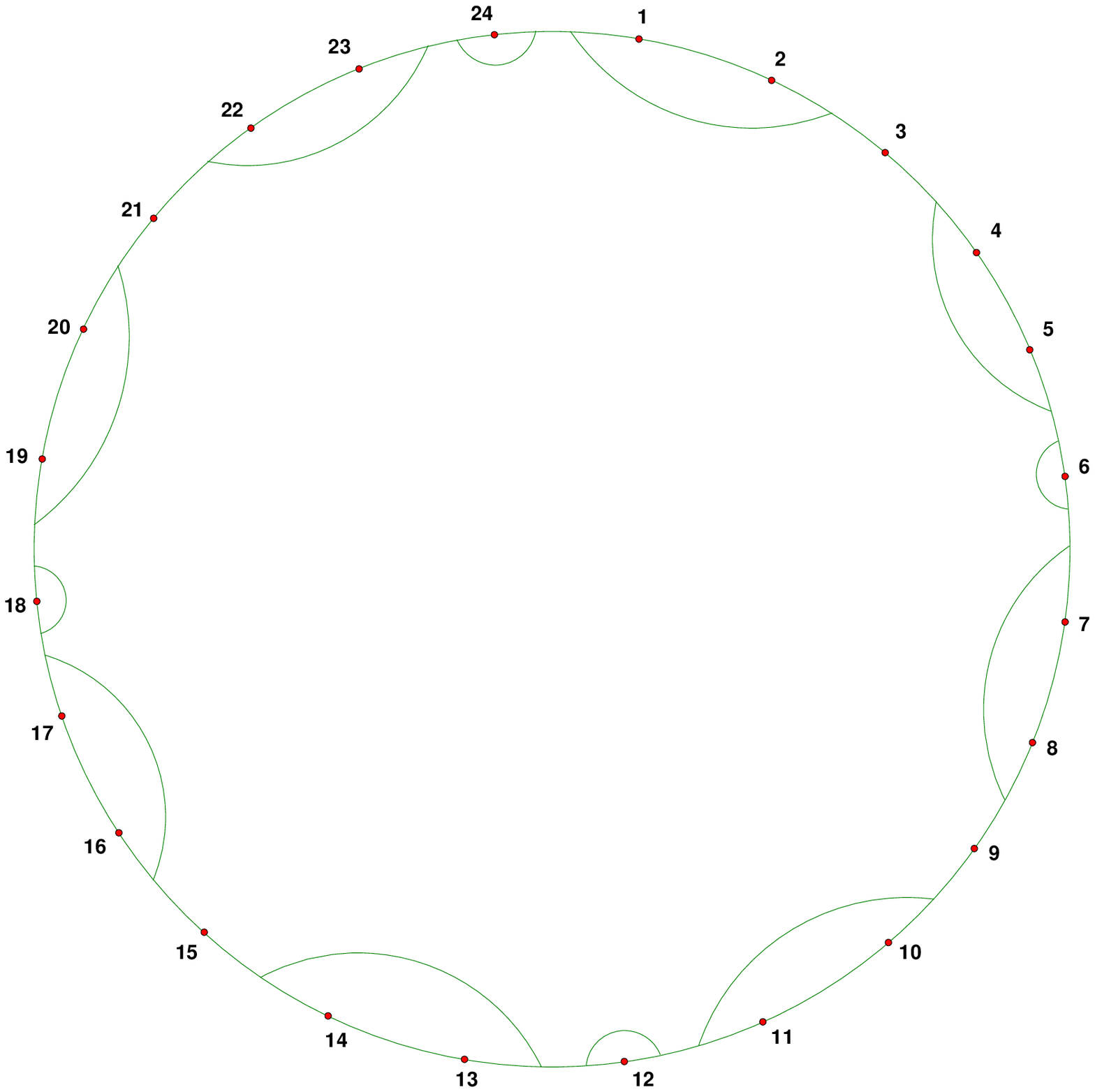}
\caption{A 6-clickable partition, not 3-clickable}
\label{F4}
\end{minipage}
\hfill
\end{figure}

Figure \ref{F3} displays a partition of $[24]$ which is $3$-, $6$-, and $12$-clickable (we use partitioning walls instead of polygons for viewability). The generator (click) $\sigma$ of $\mathbb{Z}_{24}$ may be thought of as a rotation of the diagram through $2\pi /24$ leaving the labels in place. The pictured partition is then a fixed point of $\sigma^{3i}$, $i=0,\dots ,7$. As one of the $\binom{6}{3}$ $3$-clickables, it is a fixed point of all of these, including $\sigma^3$, $\sigma^9$, $\sigma^{15}$, and $\sigma^{21}$. As one of the superset of $\binom{12}{6}$ $6$-clickables, (cf Figure \ref{F4}), it may not be invariant under those $4$ rotations, but it must be a fixed point of $\sigma^6$, and $\sigma^{18}$, while as one of the $\binom{24}{12}$ $12$-clickables, it need only be an invariant of $\{\sigma^{12},\sigma^0\}$. These considerations generalize succinctly in the following enumeration.

\vfill
\eject

\begin{theorem}
The number of NC Partition Patterns of $n$ points on a circle is
$$\frac{1}{n}\Biggl( \frac{1}{n+1}\binom{2n}{n}+\sum_{\substack{1\leq i<n\\ i|n}}\phi(\frac{n}{i})\binom{2i}{i}\Biggr)$$.
\label{ReinerCount}
\end{theorem}

\noindent \textbf{Proof} The Cauchy-Frobenius principle for the rotation group of order $n$ counts our equivalence classes by summing over all group elements the number of objects (labelled NC partitions) invariant under the element, then dividing by the group order. The identity element accounts for the Catalan number as the left summand. It is easily verified that those non-identity group elements which fix all $i$-clickable partitions but not all $j$-clickable partitions for $j>i$ are exactly those of order $\frac{n}{i}$, and these number $\phi(\frac{n}{i})$. This gives the right summand, using the result of V. Reiner. \qed

As shown in Section 4, the central binomial coefficient counts those NC Partition Patterns invariant by any reflection across a diameter which fixes the $n$ points. Cauchy-Frobenius allows us to count chirally inequivalent patterns by adding half the non-invariant patterns to the invariant ones.   

The enumeration of bicolored plane trees predates, and thus confirms, Theorem \ref{ReinerCount}.
The number of (free) plane trees on $n$ edges is known to be
$$\FPT(n):=\frac{1}{2n}\sum_{d \vert n}\phi(\frac{n}{d})\binom{2d}{d} - 
\frac{(C_{n}-C_{\frac{n-1}{2}})}{2},$$ 
\noindent see (\htmladdnormallink{A002995}{http://www.research.att.com:80/cgi-bin/access.cgi/as/njas/sequences/eisA.cgi?Anum=A002995}).
Here and below $C_{m}$ is understood to be 0 if $m$ is not an integer, 
and the term 
involving Catalan numbers is an integer because $C_{m}$ is 
odd if{f} the integer $m$ has the form $2^{k}-1$. 

The \emph{size} of a plane tree is its number of edges. The
\emph{subtrees} of a vertex are the plane trees obtained by deleting
the vertex and its incident edges.
A \emph{center} of a plane tree is a vertex $v$ that minimizes
$\max\{\textrm{size}(T)\,:\textrm{\,$T$ a subtree of $v$}\}$. A plane
tree either has a unique center or two adjacent centers. Deleting the 
connecting edge in
the latter case leaves two ordered trees. Symmetry then implies that the number of bicolored plane tree on $n$
edges is $2\,\FPT(n)-C_{\frac{n-1}{2}}$ 
(\htmladdnormallink{A054357}{http://www.research.att.com:80/cgi-bin/access.cgi/as/njas/sequences/eisA.cgi?Anum=A054357}), 
since ordered trees are yet another manifestation of the Catalan numbers. 
Clearly this agrees with the Cauchy-Frobenius count of NC partition patterns above. 

Michel Bousquet applied Cauchy-Frobenius (``Lemme de Burnside'') to enumerate $m$-ary cacti in \cite{BousquetThese}, applying a scheme due to Liskovets. His result includes bicolored plane trees as a special case.

Note added: A result equivalent to Theorem 3 and its consequences appeared in \cite{potts}. 
In particular, the formula we give for NC Dihedral Classes is Corollary 2.1 in that paper, 
which is, we believe, its first occurrence in the literature. 
We thank the authors of \cite{potts} for notifying us of these facts.

\begin{tabular}{ | c | c | c | c | c | }\hline
  $n$ & NC Rotation Classes & NC Dihedral Classes & NC Chiral patterns $\div 2$ \\
  $ $ &  \emph{or} NCPP's    & (chiral equivalence)&                             \\  
  $ $ & A054357             &(A054357+A001405)/2    &(A054357-A001405)/2         \\ \hline
  $1$ &1 &1 &0 \\
  $2$ &2 &2 &0 \\
  $3$ &3 &3 &0 \\
  $4$ &6 &6 &0 \\
  $5$ &10 &10 &0 \\
  $6$ &28 &24 &4 \\
  $7$ &63 &49 &14 \\
  $8$ &190 &130 &60 \\
  $9$ &546 &336 &210 \\
  $10$ &1708 &980 &728 \\
  $11$ &5346 &2904 &2442 \\
  $12$ &17428 &9176 &8252 \\
  $13$ &57148 &29432 &27716 \\
  $14$ &191280 &97356 &93924 \\
  $15$ &646363 &326399 &319964 \\
  $16$ &2210670 &1111770 &1098900 \\
  $17$ &7626166 &3825238 &3800928 \\
  $18$ &26538292 &13293456 &13244836 \\
  $19$ &93013854 &46553116 &46460738 \\
  $20$ &328215300 &164200028 &164015272 \\
  $21$ &1165060668 &582706692 &582353976 \\
  $22$ &4158330416 &2079517924 &2078812492 \\ \hline
\end{tabular}


\begin{thebibliography}{99}
    
\bibitem{Adrianov} N. Adrianov and A.ÊZvonkin, Composition of plane trees, 
\emph{Acta Applicandae Mathematicae} \textbf{52},
Numbers 1-3, July 1998, 239--245. 
    


\bibitem{BousquetCacti} Miklos Bona, Michel Bousquet, Gilbert Labelle and Pierre Leroux, 
Enumeration of m-ary cacti, \emph{Advances in Applied Mathematics}, \textbf{24} (2000), 22--56. 


\bibitem{BousquetThese} Michel Bousquet, 
Quelques r\`{e}sultats sur les cactus planaires, 
\emph{Annales des Sciences Mathematiques du Quebec}, \textbf{24} (2000) No. 2. p. 107-128. 

\bibitem{Bordeaux} Mireille Bousquet-M\'{e}lou,  Combinatorics in Bordeaux,
\htmladdnormallink{www.mat.univie.ac.at/$\,\widetilde{\ }\,$slc/wpapers/s34bordeaux\_des.PS}{www.mat.univie.ac.at/~slc/wpapers/s34bordeaux_des.PS}.



\bibitem{conjpart} David Callan, On conjugates for integer compositions 
and set partitions, preprint, \htmladdnormallink{http://front.math.ucdavis.edu/math.CO/0508052}{http://front.math.ucdavis.edu/math.CO/0508052}
.    

\bibitem{potts} Shu-Ciuan Chang, Jesper Lykke Jacobsen, Jesus Salas and Robert Shrock,
Exact Potts Model Partition Functions for Strips of the Triangular Lattice, 
\emph{J. Statistical Physics}, \textbf{114} (2004),  763--823. 


\bibitem{Deutsch} Emeric Deutsch, Dyck Path Enumeration,
\emph{Discrete Math.} \textbf{204} (1999),  167--202.

    
\bibitem{krewerasNC} Germain Kreweras, Sur les partitions non croisees d'un cycle,
\emph{Discrete Math.} \textbf{1}  (1972), no. 4, 333--350.


\bibitem{MotzkinBAMS} Theodore Samuel Motzkin,
Relations between hypersurface cross ratios, and a
 combinatorial formula for partitions of a polygon, for permanent
 preponderance, and for non-associative products, 
\emph{Bull. Amer. Math. Soc.}, \textbf{54} (1948), 352--360.


\bibitem{reinerNC} Victor Reiner, 
Noncrossing partitions for classical rotation groups,
\emph{Discrete Math.} \textbf{177} (1997),  195--222.


\bibitem{simionNC} Rodica Simion and  Daniel Ullman, 
On the structure of the lattice of noncrossing partitions,
\emph{Discrete Math.} \textbf{98} (1991), no. 3, 193Ð-206.

    
\bibitem{ec2} Richard P.~Stanley, \emph{Enumerative Combinatorics 
Vol.\,2}, Cambridge University Press, 1999. Exercise 6.19 and related 
material on Catalan numbers are available 
online at
\htmladdnormallink{http://www-math.mit.edu/$\,\widetilde{\ }\,$rstan/ec/ }{http://www-math.mit.edu/~rstan/ec/}.


\end{thebibliography}
\end{document}